\documentclass[reqno]{article}
\usepackage{microtype}
\usepackage{amsmath,amsfonts,amsthm,amssymb}
\usepackage{bm}
\usepackage{indentfirst}
\usepackage{graphicx}
\numberwithin{equation}{section}
\usepackage[margin=2.9cm]{geometry}
\usepackage[colorlinks,linkcolor=blue,citecolor=blue]{hyperref}
\usepackage{xcolor}

\theoremstyle{plain}
\newtheorem{Thm}{Theorem}[section]
\newtheorem*{Thm*}{Theorem}
\newtheorem{Lem}[Thm]{Lemma}

\newtheorem{Prop}[Thm]{Proposition}

\theoremstyle{definition}

\newtheorem{Rem}[Thm]{Remark}
\newtheorem{Ex}[Thm]{Example}

\newcommand{\p}{\partial}

\newcommand{\R}{\mathbb{R}}
\newcommand{\I}{\mathrm{I}}

\begin{document}

\title{Asymptotic Stability and $L^p$ Time-Decay Rates of  Rarefaction Waves for Systems of Viscous Conservation Laws}

\author{Mutong He\thanks{Academy of Mathematics and Systems Science, Chinese Academy of Sciences, Beijing 100190, PR China. E-mail: hemutong25@mails.ucas.ac.cn},
Feimin Huang\thanks{Academy of Mathematics and Systems Science, Chinese Academy of Sciences, Beijing 100190, PR China. E-mail: fhuang@amt.ac.cn},
Lingjun Liu\thanks{School of Mathematics, Statistics and Mechanics, Beijing University of Technology, Beijing 100124, PR China. E-mail: lingjunliu@bjut.edu.cn},
Lingda Xu\thanks{Hetao Institute of Mathematics and Interdisciplinary Sciences, Shenzhen 518048, PR China. E-mail: xuldmath@gmail.com}
}

\date{}
\maketitle

\begin{abstract}
We study the asymptotic stability and quantitative time-decay rates of rarefaction waves for one-dimensional $n\times n$ strictly hyperbolic systems of viscous conservation laws. While the nonlinear stability of rarefaction waves for $2\times 2$ systems was established by Xin \cite{Xin1988,Xin1989}, extending such results to general $n\times n$ systems presents additional difficulties arising from the possible coexistence of positive and negative transverse relative characteristic speeds and nonlinear interactions between distinct characteristic families. Moreover, establishing standalone $L^p$ energy estimates and explicit algebraic decay rates for general viscous systems without $L^1$ localization remains a significant challenge.

In this paper, we establish both the global $H^2(\R)$ asymptotic stability and, for the first time, standalone physical-space $L^p$ energy estimates and explicit algebraic $L^p$ time-decay rates for weak rarefaction waves. Specifically, for every integer $p\ge 3$, we prove that the perturbation $U(t)-\bar{U}(t)$ around the smooth Burgers rarefaction profile $\bar{U}$ decays in $L^p(\R)$ at the rate
\[
\|U(t)-\bar{U}(t)\|_{L^p}\le C_p\bigl(\|U_0-\bar{U}(0)\|_{H^2}^2+\delta^3\bigr)^{1/2}(1+t)^{-\frac{1}{2}(\frac{1}{2}-\frac{1}{p})},
\]
where $\delta$ denotes the wave strength. Our proof combines a diagonal characteristic weight matrix that yields coercive rarefaction damping in all characteristic modes with directed wave-interaction functionals that exploit the strict separation of characteristic speeds to control nonlinear coupling between distinct characteristic families. These estimates lead to a closed $L^p$ Lyapunov inequality, from which the algebraic time-decay rates follow by the one-dimensional Nash inequality.
\end{abstract}



\tableofcontents

\section{Introduction}\label{sec1}
We consider the Cauchy problem for a one-dimensional system of viscous conservation laws
\begin{align}\label{eq1.1}
U_t+F(U)_x=\nu U_{xx},\qquad (t,x)\in(0,\infty)\times\R,\qquad U(t,x)\in\R^n,
\end{align}
where $F:\mathcal{O}\subset\R^n\to\R^n$ is a smooth strictly hyperbolic flux and $\nu>0$ is the viscosity coefficient. By the scaling $(t,x)=(\nu\tau,\nu y)$, \eqref{eq1.1} reduces to the normalized system
\begin{align}\label{eq1.2}
U_t+F(U)_x=U_{xx}.
\end{align}
Our goal is to investigate the global asymptotic stability and explicit $L^p(\R)$ time-decay rates of a weak rarefaction wave associated with a genuinely nonlinear interior characteristic family $2\le k\le n-1$.

The large-time stability of rarefaction waves for viscous conservation laws and compressible fluid models has been a central subject since the pioneering works of Il'in--Oleinik \cite{IlinOleinik1960} for scalar equations and Matsumura--Nishihara \cite{MatsumuraNishihara1986,MatsumuraNishihara1992} and Liu--Xin \cite{LiuXin1988} for the compressible Navier--Stokes equations. For general $2\times 2$ viscous conservation laws, Xin \cite{Xin1988,Xin1989} established the nonlinear $H^2(\R)$ stability of weak rarefaction waves; see also \cite{KawashimaMatsumuraNishihara1986,NishiharaYangZhao2004,YangZhao2005,HuangLiWang2012}. However, in any $2\times 2$ strictly hyperbolic system, both characteristic families are extremal ($k=1$ or $k=2$), so all transverse characteristic speeds lie on one side of the active rarefaction speed. For general $n\times n$ systems ($n\ge 3$) with an interior rarefaction wave ($2\le k\le n-1$), or for obtaining quantitative time-decay rates without $L^1$ localization, two fundamental analytical difficulties arise:

First, at the linear level, the transverse characteristic speeds lie on both sides of the active rarefaction speed $\lambda_k$:
\[
\lambda_i-\lambda_k<0 \quad (i<k),\qquad \lambda_i-\lambda_k>0 \quad (i>k).
\]
In general, a single scalar weight is insufficient to guarantee positive rarefaction damping in all transverse characteristic modes. For instance, consider the decoupled \(3\times 3\) model
\[
F(u,v,r)=\Bigl(-u,\,\frac{1}{2}v^2,\,r\Bigr)^T,\qquad \bar{U}=(0,w,0)^T,\qquad |w|<1,
\]
testing with a scalar weight $d(w)$ yields the transverse drift coefficients $(w+1)d'(w)$ and $(w-1)d'(w)$, which have opposite signs on $|w|<1$. Moreover, as shown in Examples~\ref{ex3.3} and \ref{ex3.4} below, neither does the existence of a strictly convex entropy guarantee positive transverse damping, nor does a general strictly hyperbolic flux necessarily admit a local $C^3$ strictly convex entropy.

Second, at the nonlinear level, writing the perturbation in characteristic coordinates $W=U-\bar{U}=R(w)Z$ produces quadratic flux integrals of the form
\[
\int_{\R} T_{abc}(w)z_{a,x}z_b z_c\,dx.
\]
When $a=b=c$, integration by parts transfers the spatial derivative onto the background coefficient and gains a factor of $w_x>0$, which is absorbed by the rarefaction damping. When two distinct characteristic families appear in $(a,b,c)$, however, the instantaneous $L^2$ or $L^p$ componentwise energy provides no cross-family dissipation to control the product. For instance, in the $3\times 3$ system
\[
F(u,v,r)=\Bigl(-u+\frac{\mu}{2}r^2,\,\frac{1}{2}v^2,\,r\Bigr)^T,\qquad \mu\neq 0,
\]
around $\bar{U}=(0,w,0)^T$, the transverse perturbations with $v=0$ satisfy $u_t-u_x+\frac{\mu}{2}(r^2)_x=u_{xx}$ and $r_t+r_x=r_{xx}$. Multiplying by $|u|^{p-2}u$ and $|r|^{p-2}r$ gives
\[
\frac{1}{p}\frac{d}{dt}\bigl(\|u\|_{L^p}^p+\|r\|_{L^p}^p\bigr)+(p-1)\int_{\R}\bigl(|u|^{p-2}u_x^2+|r|^{p-2}r_x^2\bigr)\,dx=\frac{\mu(p-1)}{2}\int_{\R}|u|^{p-2}u_x r^2\,dx.
\]
Under the $L^2$-conserving scaling $(u_\ell,r_\ell)(x)=\varepsilon\ell^{-1/2}(\phi,\psi)(x/\ell)$, the right-hand side of this identity scales as $\varepsilon^{p+1}\ell^{-(p+1)/2}$, whereas the left-hand side dissipation scales as $\varepsilon^p\ell^{-(p+2)/2}$; for $\ell\gg\varepsilon^{-2}$ and waves translated away from the rarefaction fan, the cubic flux term cannot be absorbed by the instantaneous diffusion or fan damping. This explains why standalone $L^p$ ($p>2$) energy estimates---which are classical for scalar viscous conservation laws \cite{HattoriNishihara1991,Ito1996,KawashimaNishibataNishikawa2004,HuangXu2022}---were previously unavailable for systems around rarefaction waves.

In this paper, we resolve both obstacles and establish, for the first time, the global $H^2(\R)$ stability of weak interior rarefaction waves together with standalone physical-space $L^p$ ($p\ge 4$) energy estimates and explicit algebraic time-decay rates for $p\ge 3$. Under the initial perturbation bound $A:=(\|U_0-\bar{U}(0)\|_{H^2}^2+\delta^3)^{1/2}\ll 1$, we prove in Theorem~\ref{thm2.1} that \eqref{eq1.2} admits a unique global solution satisfying
\[
\|U(t)-\bar{U}(t)\|_{L^p}\le C_p A(1+t)^{-\frac{1}{2}(\frac{1}{2}-\frac{1}{p})}\qquad \forall\,t\ge 0,\quad p\ge 3,
\]
without any $L^1$ localization, zero-mass assumption, or anti-derivative transformation.

Our approach relies on three main ingredients in physical space:
\begin{enumerate}
\item \textit{Diagonal characteristic weight matrix (Section~\ref{sec3}).} In the characteristic variables $Z=L(w)(U-\bar{U})$, we construct a diagonal matrix weight $D(w)=\operatorname{diag}(d_1(w),\dots,d_n(w))$ whose entries below and above the active $k$-th family satisfy scalar ODEs tilting in opposite directions. This makes the transverse drift matrix positive definite in every characteristic family, while genuine nonlinearity provides the active $k$-th component damping.
\item \textit{Directed quartic wave-interaction functional for $H^2$ closure (Sections~\ref{sec4}--\ref{sec5}).} To control the cross-family cubic flux integrals $\int_\R T_{abc}(w)z_{a,x}z_b z_c\,dx$ without $L^1$ assumptions, we introduce for $i<j$ the half-space interaction functional
\[
\mathcal{I}_{ij}(t)=\iint_{x>y} z_i(t,x)^2 z_j(t,y)^2\,dx\,dy,
\]
which is a continuous analogue of the Glimm--Bony interaction potential \cite{Glimm1965,Bony1987}. By strict hyperbolicity ($\lambda_j-\lambda_i\ge\theta_0>0$), differentiating $\mathcal{I}_{ij}(t)$ in time produces the coercive cross-family spacetime integral
\[
\int_0^T\int_{\R} z_i(t,x)^2 z_j(t,x)^2\,dx\,dt,
\]
which absorbs all non-resonant cubic flux terms and closes the global $H^2(\R)$ energy estimate.
\item \textit{Standalone $L^p$ energy and $(p,2)+(2,p)$ wave interactions (Sections~\ref{sec6}--\ref{sec9}).} For each fixed integer $p\ge 4$, we construct a $p$-adapted characteristic weight $D(w)$ and couple the weighted $L^p$ energy $\mathcal{E}_p(t)$ with both directed mixed-power interaction functionals
\[
\mathcal{I}_{ij}^{p,2}(t)=\iint_{x>y}|z_i(t,x)|^p z_j(t,y)^2\,dx\,dy,\qquad \mathcal{I}_{ij}^{2,p}(t)=\iint_{x>y} z_i(t,x)^2|z_j(t,y)|^p\,dx\,dy.
\]
Using the algebraic bound $a^{p-2}b^4\le a^p b^2+a^2 b^p$ ($p\ge 4$) and integrating the repeated-factor nonlinear flux terms by parts prior to taking $L^\infty$ bounds on the cumulative tails, we obtain a closed differential inequality for the $p$-level Lyapunov functional $\mathcal{H}_p(t)=\mathcal{E}_p(t)+K_p\mathcal{I}_p(t)$. Applying the one-dimensional Nash inequality \cite{Nash1958} then converts the $p$-diffusion $\mathcal{D}_p(t)=\sum_{i=1}^n\int_\R|z_i|^{p-2}z_{i,x}^2\,dx$ into a nonlinear decay ODE, yielding the algebraic $L^p$ decay rate $(1+t)^{-\frac{1}{2}(\frac{1}{2}-\frac{1}{p})}$.
\end{enumerate}

Let us place our results in the context of the existing literature. For inviscid hyperbolic systems, $BV/L^1$ stability conditions for rarefaction waves were studied by Lewicka \cite{Lewicka2005a,Lewicka2005b}, and vanishing viscosity semigroups for small-$BV$ data were constructed by Bianchini and Bressan \cite{BianchiniBressan2005}. For general viscous conservation laws, Szepessy--Zumbrun \cite{SzepessyZumbrun1996} proved stability of weak rarefaction waves under localized perturbations via pointwise Green's function iterations, and Liu--Yu \cite{LiuYu2010,LiuYu2022} established detailed pointwise wave structures using Green's functions and inner--outer expansions; see also \cite{Zeng1994,LiuZeng1997,MasciaZumbrun2004,Zeng2018,Zeng2019} for $L^1$/pointwise decay around constant states and shock profiles, and \cite{MatsumuraNishihara1985,HuangMatsumuraXin2006,HuangXinYang2008,HuangLiMatsumura2010,Meng2025} for viscous shocks, contact waves, and composite waves. Physical-space $L^p$ energy methods for scalar viscous conservation laws were developed in \cite{HattoriNishihara1991,Ito1996,NishikawaNishihara2000,KawashimaNishibataNishikawa2004,HuangXu2022}, and a physical-space hypocoercive energy approach for partially dissipative systems near constant states was recently introduced in \cite{Crin-BaratShouZuazua2025}. To the best of our knowledge, Theorem~\ref{thm2.1} provides the first global $H^2$ stability result for interior rarefaction waves of general viscous systems without a convex entropy, as well as the first standalone $L^p$ energy estimate and explicit time-decay rate for rarefaction waves in the $H^2(\R)$ framework.

The paper is organized as follows. Section~\ref{sec2} introduces the Burgers rarefaction profile and states the main theorem (Theorem~\ref{thm2.1}). Section~\ref{sec3} constructs the diagonal characteristic weight matrix and presents three illustrative $3\times 3$ examples. Sections~\ref{sec4} and \ref{sec5} derive the characteristic perturbation equations and close the global $H^2$ estimate via the quartic interaction functional. Sections~\ref{sec6} and \ref{sec7} establish the weighted $L^p$ energy estimate and the $(p,2)+(2,p)$ wave-interaction estimates for $p\ge 4$. Sections~\ref{sec8} and \ref{sec9} close the $L^p$ Lyapunov inequality and derive the $L^p$ decay rates via the Nash inequality. Finally, Section~\ref{sec10} completes the local existence, global continuation, and $L^\infty$ convergence to the inviscid rarefaction wave.

\section{Smooth rarefaction profile and main results}\label{sec2}
Let $F\in C^7(\mathcal{O};\R^n)$ on an open state domain $\mathcal{O}\subset\R^n$, and denote $A_F(U):=DF(U)$. On a compact tubular neighborhood $\mathcal{K}\Subset\mathcal{O}$, we assume that $A_F(U)$ is strictly hyperbolic with real eigenvalues satisfying
\begin{align}\label{eq2.1}
\lambda_1(U)<\cdots<\lambda_n(U),\qquad \min_{i\neq j}|\lambda_i(U)-\lambda_j(U)|\ge\theta_0>0,
\end{align}
and smooth dual left and right eigenvectors $l_i(U), r_j(U)$ normalized by $l_i(U)r_j(U)=\delta_{ij}$. Fix an interior family $2\le k\le n-1$ and assume that the $k$-th characteristic field is genuinely nonlinear on $\mathcal{K}$:
\[
\nabla\lambda_k(U)\cdot r_k(U)\ge g_0>0.
\]
After rescaling $r_k$ so that $\nabla\lambda_k(U)\cdot r_k(U)\equiv 1$ along the $k$-rarefaction curve $U_c(w)$, we parametrize the curve by the characteristic speed $w$:
\[
U_c'(w)=r_k(U_c(w)),\qquad \lambda_k(U_c(w))=w.
\]
Fix asymptotic end states $w_-<w_+$, set $\delta:=w_+-w_-$ and $U_\pm:=U_c(w_\pm)$, and assume that the rarefaction segment lies strictly inside $\mathcal{K}$:
\[
K_r:=\{U_c(w):w_-\le w\le w_+\}\Subset\operatorname{int}\mathcal{K},\qquad \kappa:=\operatorname{dist}(K_r,\p\mathcal{K})>0.
\]

To construct a smooth rarefaction profile with exponential spatial tails, let $w(t,x)$ solve the viscous Burgers equation
\begin{align}\label{eq2.2}
w_t+ww_x=w_{xx},\qquad w(0,x)=\frac{w_-+w_+}{2}+\frac{\delta}{2}\tanh(\delta x),\qquad w_x>0.
\end{align}
By the standard Cole--Hopf representation and comparison principle \cite{IlinOleinik1960,Xin1988,HattoriNishihara1991,LiuYu2010}, $w(t,x)$ satisfies the decay bounds
\begin{align}
\int_{\R}w_x(t,x)\,dx&=\delta,\qquad \|w_x(t)\|_{L^\infty}\le\frac{C\delta^2}{1+\delta^2 t},\label{eq2.3}\\
|\p_x^{m+1}w(t,x)|&\le C_m\delta^m w_x(t,x),\qquad 0\le m\le 4,\label{eq2.4}\\
\int_0^\infty\int_{\R}w_x^3\,dx\,dt&\le C\delta^3.\label{eq2.5}
\end{align}
More generally, for every $r\ge 2$,
\begin{align}\label{eq2.6}
\int_{\R}w_x(t,x)^r\,dx\le \delta\|w_x(t)\|_{L^\infty}^{r-1}\le C_r\delta^{2r-1}(1+\delta^2 t)^{-(r-1)},
\end{align}
and consequently,
\begin{align}\label{eq2.7}
\int_0^\infty\|w_x(t)\|_{L^\infty}^2\,dt\le C\delta^2,\quad \int_0^\infty\int_{\R}w_x^r\,dx\,dt\le C_r\delta^{2r-3}\ (r>2),\quad \int_{\R}w_x^{p+1}\,dx\le C_p\delta^{2p+1}(1+\delta^2 t)^{-p}.
\end{align}
We define the smooth $k$-rarefaction profile by
\begin{align}\label{eq2.8}
\bar{U}(t,x):=U_c(w(t,x)).
\end{align}
Since $A_F(\bar{U})U_c'(w)=wU_c'(w)$ and $w_t+ww_x=w_{xx}$, the profile residual is given explicitly by
\[
\bar{U}_t+F(\bar{U})_x-\bar{U}_{xx}=-U_c''(w)w_x^2.
\]

Our main result is stated as follows.

\begin{Thm}\label{thm2.1}
There exist constants $0<\delta_0\le 1$ and $\varepsilon_0>0$, depending only on $F$ and its eigenframe on $\mathcal{K}$, $\theta_0, g_0$, and $\kappa$, such that if $0<\delta\le\delta_0$ and
\[
W_0:=U_0-\bar{U}(0,\cdot)\in H^2(\R),\qquad \|W_0\|_{H^2}\le\varepsilon_0,
\]
then the Cauchy problem \eqref{eq1.2} with $U(0,x)=U_0(x)$ admits a unique global strong solution satisfying
\[
W:=U-\bar{U}\in C([0,\infty);H^2(\R)),\qquad W_x\in L^2(0,\infty;H^2(\R)),
\]
and the uniform $H^2$ energy estimate
\begin{align}\label{eq2.9}
\sup_{t\ge 0}\|W(t)\|_{H^2}^2+c\int_0^\infty\Bigl(\|W_x(t)\|_{H^2}^2+\sum_{m=0}^2\int_{\R}w_x|\p_x^m W|^2\,dx\Bigr)\,dt\le C\bigl(\|W_0\|_{H^2}^2+\delta^3\bigr).
\end{align}
Furthermore, set
\[
A:=\bigl(\|W_0\|_{H^2}^2+\delta^3\bigr)^{1/2}.
\]
For every fixed integer $p\ge 4$, there exist thresholds $0<\delta_p\le\delta_0$ and $0<\varepsilon_p\le\varepsilon_0$ such that whenever $\delta\le\delta_p$ and $\|W_0\|_{H^2}\le\varepsilon_p$, the solution satisfies the global $L^p$ Lyapunov estimate
\begin{align}\label{eq2.10}
\sup_{t\ge 0}\|W(t)\|_{L^p}^p+\int_0^\infty\bigl(\mathcal{D}_p(t)+\mathcal{G}_p(t)+\mathcal{Y}_p(t)+\mathcal{B}_p(t)\bigr)\,dt\le C_p A^p,
\end{align}
where $\mathcal{D}_p, \mathcal{G}_p, \mathcal{Y}_p$ are defined in \eqref{eq6.3} and $\mathcal{B}_p$ in \eqref{eq7.2}, and the algebraic $L^p$ time-decay estimate
\begin{align}\label{eq2.11}
\|W(t)\|_{L^p}\le C_p A(1+t)^{-\beta_p},\qquad \beta_p=\frac{1}{2}\Bigl(\frac{1}{2}-\frac{1}{p}\Bigr)\qquad \forall\,t\ge 0.
\end{align}
In addition, for $\delta\le\delta_4$ and $\|W_0\|_{H^2}\le\varepsilon_4$, interpolation between $L^2(\R)$ and $L^4(\R)$ yields $\|W(t)\|_{L^3}\le CA(1+t)^{-1/12}$, and if $U^R(x/t)$ denotes the centered inviscid $k$-rarefaction wave connecting $U_-$ to $U_+$, then
\begin{align}\label{eq2.12}
\lim_{t\to\infty}\sup_{x\in\R}|U(t,x)-U^R(x/t)|=0.
\end{align}
\end{Thm}

\begin{Rem}[On the characteristic family]
Theorem \ref{thm2.1} applies to rarefaction waves associated with any genuinely nonlinear characteristic family \(1\le k\le n\). Indeed, the construction of the characteristic weights in Section 3 requires only the strict separation (\ref{eq2.1}), while the directed wave-interaction estimates rely solely on the strict ordering of distinct characteristic speeds. Neither argument requires the active family to be interior. We emphasize the case \(2\le k\le n-1\) because the transverse characteristic speeds then lie on both sides of the rarefaction speed, creating a two-sided sign difficulty absent in the extremal cases \(k=1\) and \(k=n\).
\end{Rem}

\begin{Rem}
The term $\delta^3$ in $A^2=\|W_0\|_{H^2}^2+\delta^3$ arises from the spacetime $L^2$ norm of the Burgers profile residual $-U_c''(w)w_x^2$ via \eqref{eq2.5}. For $p=3$, the direct $(p,2)+(2,p)$ algebraic reduction in \eqref{eq6.16} requires $p\ge 4$; hence the $L^3$ decay rate $(1+t)^{-1/12}$ is obtained by interpolating the uniform $L^2$ bound \eqref{eq2.9} with the $p=4$ decay estimate $\|W(t)\|_{L^4}\le CA(1+t)^{-1/8}$. For the unscaled system \eqref{eq1.1} with viscosity $\nu>0$, the time factor $(1+t)$ in \eqref{eq2.11} becomes $(1+t/\nu)$.
\end{Rem}

\section{Characteristic weight matrix and model examples}\label{sec3}
Along the rarefaction curve $U_c(w)$, we denote the right and left eigenframe matrices and the eigenvalue matrix by
\[
R(w)=(r_1,\dots,r_n)(U_c(w)),\quad L(w)=R(w)^{-1},\quad \Lambda(w)=L(w)A_F(U_c(w))R(w)=\operatorname{diag}(\lambda_1,\dots,\lambda_n),
\]
and set
\[
\Omega(w):=L(w)R'(w),\qquad C(w):=\Omega(w)\bigl(\Lambda(w)-w\I\bigr)+\Lambda'(w),
\]
where primes denote derivatives with respect to $w$. For $w\in[w_-,w_+]$ and diagonal matrices $\frac{1}{2}\I\le D\le 2\I$, we define the structural constants
\[
L_*:=\sup_{1\le i\le n,\,w}|\lambda_i'(w)|,\qquad K_*:=\sup_{D,\,w}|\operatorname{off}(DC+C^T D)|,
\]
where $\operatorname{off}(M)$ denotes the off-diagonal part of a matrix $M$ and $|\cdot|$ denotes the matrix operator norm.

\begin{Lem}\label{lem3.1}
Choose
\begin{align}\label{eq3.1}
\Gamma\ge 16\bigl(1+L_*+K_*+K_*^2\bigr)
\end{align}
and take $\delta_0>0$ small enough that
\begin{align}\label{eq3.2}
\delta_0\sup_{i\neq k,\,w}\Bigl|2\Omega_{ii}(w)+\frac{\lambda_i'(w)-\Gamma}{\lambda_i(w)-w}\Bigr|\le\log 2.
\end{align}
Then for any $0<\delta\le\delta_0$, there exists a smooth diagonal matrix $D(w)=\operatorname{diag}(d_1(w),\dots,d_n(w))$ satisfying
\begin{align}\label{eq3.3}
\frac{1}{2}\I\le D(w)\le 2\I,\qquad w\in[w_-,w_+],
\end{align}
such that for $m=0,1,2$, the symmetric matrix
\begin{align}\label{eq3.4}
P_m(w):=DC+C^T D+wD'-(D\Lambda)'+2mD\Lambda'\ge c_D\I
\end{align}
is uniformly positive definite on $[w_-,w_+]$ for a constant $c_D>0$.
\end{Lem}

\begin{proof}
Set $d_k(w)\equiv 1$ and, for each $i\neq k$, solve the linear ODE
\begin{align}\label{eq3.5}
\frac{d_i'(w)}{d_i(w)}=2\Omega_{ii}(w)+\frac{\lambda_i'(w)-\Gamma}{\lambda_i(w)-w},\qquad d_i(w_-)=1.
\end{align}
Since $|\lambda_i(w)-w|=|\lambda_i-\lambda_k|\ge\theta_0>0$ by \eqref{eq2.1}, the condition \eqref{eq3.2} immediately implies \eqref{eq3.3}. A direct calculation using $C_{ii}=\Omega_{ii}(\lambda_i-w)+\lambda_i'$, $\lambda_k(w)=w$, and $\lambda_k'(w)=1$ gives the entries of $P_0$:
\[
(P_0)_{kk}=1,\qquad (P_0)_{ii}=\Gamma d_i\quad (i\neq k),\qquad (P_0)_{ij}=d_i\Omega_{ij}(\lambda_j-w)+d_j\Omega_{ji}(\lambda_i-w)\quad (i\neq j).
\]
For $m\in\{0,1,2\}$, the active $k$-th diagonal entry of $P_m$ is $1+2m\ge 1$, while each transverse diagonal entry ($i\neq k$) is $d_i(\Gamma+2m\lambda_i')\ge\frac{1}{2}(\Gamma-4L_*)$. By \eqref{eq3.1}, after including the transverse off-diagonal entries, the $(n-1)\times(n-1)$ transverse block $T$ of $P_m$ satisfies $T\ge\frac{\Gamma}{4}\I$, while the coupling vector $b$ between the $k$-th mode and the transverse block satisfies $|b|\le K_*$. Writing $P_m$ in block form with respect to the active scalar component $x\in\R$ and transverse vector $y\in\R^{n-1}$, we obtain by Cauchy--Schwarz and \eqref{eq3.1} that
\[
\begin{pmatrix}x & y^T\end{pmatrix}\begin{pmatrix}1+2m & b^T\\ b & T\end{pmatrix}\begin{pmatrix}x\\ y\end{pmatrix}\ge \frac{1}{2}x^2+\Bigl(\frac{\Gamma}{4}-2K_*^2\Bigr)|y|^2\ge c_D\bigl(x^2+|y|^2\bigr),
\]
which proves \eqref{eq3.4}.
\end{proof}

We next present three explicit $3\times 3$ examples illustrating the role of the diagonal weight matrix $D(w)$ and the independence of Theorem~\ref{thm2.1} from convex entropy assumptions; related characteristic weight techniques in $L^1/BV$ and shock stability settings can be found in \cite{Liu1985,Lewicka2005a,Lewicka2005b}.

\begin{Ex}\label{ex3.2}
Consider the decoupled $3\times 3$ system with flux $F(u,v,r)=(-u,\frac{1}{2}v^2,r)^T$ around the $2$-rarefaction profile $\bar{U}=(0,w,0)^T$ with $|w|<1$. Here $\Omega=0$, $\lambda_1=-1$, $\lambda_2=w$, and $\lambda_3=1$. If one uses a scalar weight $D(w)=d(w)\I$, the two transverse diagonal entries of $P_0$ are
\[
(P_0)_{11}=(w+1)d'(w),\qquad (P_0)_{33}=(w-1)d'(w).
\]
Since $w+1>0$ and $w-1<0$ on $|w|<1$, these two coefficients cannot both be positive for any scalar function $d(w)$. In contrast, the ODE \eqref{eq3.5} chooses $d_1'(w)>0$ and $d_3'(w)<0$, making both transverse entries strictly positive.
\end{Ex}

\begin{Ex}\label{ex3.3}
Consider the symmetric gradient flux $F=\nabla\Phi$ with
\[
\Phi(u,v,r)=\frac{1}{6}u^3+\frac{a_0}{2}v^2-\frac{c_0}{2}uv^2+\frac{\beta}{6}v^3+\frac{b_0}{2}r^2+\frac{\gamma}{6}r^3,
\]
where $a_0<0<b_0$ and $c_0>0$. Since $DF=D^2\Phi$ is symmetric, $\eta(U)=\frac{1}{2}|U|^2$ is a global strictly convex entropy for \eqref{eq1.2}. Along the interior $u$-rarefaction profile $(u,0,0)^T$ ($0\le u\le\delta\ll 1$), the linearized transverse $v$-equation is $v_t+(a_0-c_0 u)v_x-c_0 u_x v=v_{xx}$, so the unweighted $L^2$ energy identity gives
\[
\frac{1}{2}\frac{d}{dt}\|v\|_{L^2}^2+\|v_x\|_{L^2}^2-\frac{c_0}{2}\int_{\R}u_x v^2\,dx=0.
\]
For $c_0>0$ and $u_x>0$, the rarefaction drift term has the unfavorable sign, showing that even a strictly convex entropy with symmetric $DF$ does not automatically provide coercive transverse damping without the characteristic weight $D(w)$.
\end{Ex}

\begin{Ex}\label{ex3.4}
Consider $U=(u,v,r)^T$ and the triangular flux
\[
F(U)=\Bigl(-u,\,\frac{1}{2}v^2,\,r+uv\Bigr)^T,\qquad DF(U)=\begin{pmatrix}-1 & 0 & 0\\ 0 & v & 0\\ v & u & 1\end{pmatrix}.
\]
On $|v|<1$, the eigenvalues are $-1<v<1$, the middle family is genuinely nonlinear with $\lambda_2(U)=v$, $r_2(U)=(0,1,-\frac{u}{1-v})^T$, $\nabla\lambda_2\cdot r_2=1$, and $U_c(w)=(0,w,0)^T$ is a $2$-rarefaction curve. However, this system admits no $C^3$ strictly convex entropy $\eta$ in any neighborhood of the origin. Indeed, if $H=D^2\eta$ symmetrized $DF$, then since $DF(0)=\operatorname{diag}(-1,0,1)$ has distinct eigenvalues, $H(0)$ would be diagonal with $h_3:=H_{33}(0)>0$. The $(1,2)$ symmetry entry of $HDF$ reads $(1+v)H_{12}+uH_{13}-vH_{23}=0$; differentiating with respect to $r$ at the origin gives $\p_{uvr}^3\eta(0)=0$. On the other hand, the $(1,3)$ symmetry entry reads $2H_{13}=vH_{33}$, whose $v$-derivative at the origin gives $2\p_{uvr}^3\eta(0)=h_3>0$, a contradiction. Theorem~\ref{thm2.1} nevertheless applies directly to this system.
\end{Ex}

\section{Weighted \texorpdfstring{$H^2$}{H2} energy estimates}\label{sec4}
In terms of the characteristic perturbation $Z=L(w)(U-\bar{U})=(z_1,\dots,z_n)^T$, we denote $Z_m:=\p_x^m Z$ for $m=0,1,2$ and define the weighted $H^2$ energy and dissipation functionals
\begin{align*}
\mathcal{E}_2(t)&:=\frac{1}{2}\sum_{m=0}^2\int_{\R}Z_m^T D(w)Z_m\,dx,\\
\mathcal{D}_2(t)&:=\sum_{m=0}^2\|Z_{m+1}\|_{L^2}^2=\|Z_x\|_{H^2}^2,\\
\mathcal{G}_2(t)&:=\sum_{m=0}^2\int_{\R}w_x|Z_m|^2\,dx.
\end{align*}

\begin{Lem}\label{lem4.1}
Assume that on $[0,T]$,
\begin{align}\label{eq4.1}
\sup_{0\le t\le T}\|Z(t)\|_{H^2}\le\varepsilon_b.
\end{align}
There exist constants $c_0,C>0$ independent of $T$ such that if $\delta+\varepsilon_b$ is sufficiently small, then
\begin{align}\label{eq4.2}
\mathcal{E}_2'(t)+c_0\bigl(\mathcal{D}_2(t)+\mathcal{G}_2(t)\bigr)\le T_3(t)+C\int_{\R}w_x(t,x)^3\,dx,
\end{align}
where $T_3(t)$ is the cubic flux integral
\begin{align}\label{eq4.3}
T_3(t):=\sum_{a,b,c=1}^n\int_{\R}T_{abc}(w)z_{a,x}z_b z_c\,dx
\end{align}
with coefficients
\begin{align}\label{eq4.4}
T_{abc}(w):=\frac{1}{2}d_a(w)l_a(U_c(w))D^2 F(U_c(w))[r_b(U_c(w)),r_c(U_c(w))].
\end{align}
Consequently, for every $0\le\tau\le T$,
\begin{align}\label{eq4.5}
\sup_{0\le t\le\tau}\mathcal{E}_2(t)+c_0\int_0^\tau\bigl(\mathcal{D}_2(t)+\mathcal{G}_2(t)\bigr)\,dt\le C\bigl(\mathcal{E}_2(0)+\delta^3\bigr)+C\sup_{0\le\sigma\le\tau}\Bigl|\int_0^\sigma T_3(t)\,dt\Bigr|.
\end{align}
\end{Lem}

\begin{proof}
Subtracting the profile equation $\bar{U}_t+F(\bar{U})_x-\bar{U}_{xx}=-U_c''(w)w_x^2$ from \eqref{eq1.2} gives the perturbation equation for $W=U-\bar{U}$:
\begin{align}\label{eq4.6}
W_t+\bigl(A_F(\bar{U})W\bigr)_x-W_{xx}=-\p_x N_0(\bar{U},W)+U_c''(w)w_x^2,
\end{align}
where $N_0(\bar{U},W):=F(\bar{U}+W)-F(\bar{U})-A_F(\bar{U})W$. Substituting $W=R(w)Z$ into \eqref{eq4.6}, multiplying on the left by $L(w)$, and using $w_t=w_{xx}-ww_x$ and $L(w)A_F(\bar{U})R(w)=\Lambda(w)$, we see that the $w_{xx}LR'Z$ terms cancel and $Z$ satisfies
\begin{align}\label{eq4.7}
Z_t+\Lambda Z_x-Z_{xx}+w_x CZ=2w_x\Omega Z_x+w_x^2\Theta Z-L\p_x N_0+w_x^2 LU_c''(w),
\end{align}
where $\Theta(w):=L(w)R''(w)$. Defining the characteristic nonlinear flux
\[
N(w,Z):=L(w)N_0(U_c(w),R(w)Z)
\]
and using $L'R=-LR'=-\Omega$, we rewrite the nonlinear term as
\begin{align}\label{eq4.8}
-L\p_x N_0=-\p_x N-w_x\Omega N.
\end{align}

Applying $\p_x^m$ ($m=0,1,2$) to \eqref{eq4.7} yields
\[
(Z_m)_t+\Lambda Z_{m+1}-Z_{m+2}+w_x CZ_m+mw_x\Lambda' Z_m=H_m,\qquad m=0,1,2,
\]
where
\begin{align}
H_0&=2w_x\Omega Z_1+w_x^2\Theta Z_0-L\p_x N_0+w_x^2 LU_c'',\label{eq4.9}\\
H_1&=2w_x\Omega Z_2+\p_x(2w_x\Omega)Z_1+w_x^2\Theta Z_1+\p_x(w_x^2\Theta)Z_0-\p_x(w_x C)Z_0-\p_x(L\p_x N_0)+\p_x(w_x^2 LU_c''),\nonumber\\
H_2&=2w_x\Omega Z_3+2\p_x(2w_x\Omega)Z_2+\p_x^2(2w_x\Omega)Z_1+w_x^2\Theta Z_2+2\p_x(w_x^2\Theta)Z_1+\p_x^2(w_x^2\Theta)Z_0\nonumber\\
&\quad-\p_x^2\Lambda\,Z_1-2\p_x(w_x C)Z_1-\p_x^2(w_x C)Z_0-\p_x^2(L\p_x N_0)+\p_x^2(w_x^2 LU_c'').\label{eq4.10}
\end{align}
Note that by \eqref{eq3.3} and \eqref{eq2.3}--\eqref{eq2.4}, $\mathcal{E}_2(t)$ and $\mathcal{D}_2(t)+\mathcal{G}_2(t)$ are uniformly equivalent to the physical perturbation norms for small $\delta$:
\[
c_E\|Z\|_{H^2}^2\le\mathcal{E}_2(t)\le C_E\|Z\|_{H^2}^2,\quad \mathcal{E}_2(t)\simeq\|W\|_{H^2}^2,\quad \mathcal{D}_2(t)+\mathcal{G}_2(t)\simeq\|W_x\|_{H^2}^2+\sum_{m=0}^2\int_{\R}w_x|\p_x^m W|^2\,dx.
\]

Taking the inner product of the $Z_m$ equation with $D(w)Z_m$, integrating by parts over $\R$, and using $D_t=D'(w)(w_{xx}-ww_x)$ gives, for each $m=0,1,2$,
\begin{align}\label{eq4.11}
&\frac{1}{2}\frac{d}{dt}\int_{\R}Z_m^T DZ_m\,dx+\int_{\R}Z_{m+1}^T DZ_{m+1}\,dx+\frac{1}{2}\int_{\R}w_x Z_m^T P_m Z_m\,dx\nonumber\\
&=\int_{\R}Z_m^T DH_m\,dx-2\int_{\R}w_x Z_m^T D'Z_{m+1}\,dx-\frac{1}{2}\int_{\R}w_x^2 Z_m^T D''Z_m\,dx.
\end{align}
By Lemma~\ref{lem3.1}, summing the left-hand side of \eqref{eq4.11} over $m=0,1,2$ provides coercive control of $\mathcal{D}_2(t)+\mathcal{G}_2(t)$. We now estimate the right-hand side of \eqref{eq4.11}.

First, by Cauchy--Schwarz and Young's inequality, for any $\eta>0$,
\begin{align}\label{eq4.12}
\Bigl|2\int_{\R}w_x Z_m^T D\Omega Z_{m+1}\,dx-2\int_{\R}w_x Z_m^T D'Z_{m+1}\,dx\Bigr|&\le\eta\|Z_{m+1}\|_{L^2}^2+C\|w_x\|_{L^\infty}\int_{\R}w_x|Z_m|^2\,dx,\nonumber\\
\Bigl|\frac{1}{2}\int_{\R}w_x^2 Z_m^T D''Z_m\,dx\Bigr|&\le C\|w_x\|_{L^\infty}\int_{\R}w_x|Z_m|^2\,dx.
\end{align}
Second, by the product rule and \eqref{eq2.3}--\eqref{eq2.4},
\begin{align*}
|\p_x^r(w_x\Omega)|+|\p_x^r(w_x C)|&\le C\delta^r w_x\quad (r=1,2),\qquad |\p_x^2\Lambda|\le C\delta w_x,\\
|\p_x^r(w_x^2\Theta)|+|\p_x^r(w_x^2 LU_c'')|&\le C\delta^r w_x^2\quad (r=0,1,2).
\end{align*}
Hence all linear background derivative terms in \eqref{eq4.9}--\eqref{eq4.10} are bounded by $C(\delta+\|w_x\|_{L^\infty})\mathcal{G}_2(t)$, while the profile residual terms satisfy
\begin{align}\label{eq4.13}
\Bigl|\int_{\R}Z_m^T D\p_x^m(w_x^2 LU_c'')\,dx\Bigr|\le C\delta^m\int_{\R}w_x^2|Z_m|\,dx\le\eta\int_{\R}w_x|Z_m|^2\,dx+C\delta^{2m}\int_{\R}w_x^3\,dx.
\end{align}

Third, using \eqref{eq4.8} and integrating the outer $\p_x$ by parts onto the multiplier, we write the nonlinear flux pairings for $m=0,1,2$ as
\begin{align}
J_0&=\int_{\R}Z_1^T DN\,dx+\int_{\R}w_x Z_0^T D'N\,dx-\int_{\R}w_x Z_0^T D\Omega N\,dx,\label{eq4.14}\\
J_1&=\int_{\R}Z_2^T D\p_x N\,dx+\int_{\R}w_x Z_1^T D'\p_x N\,dx-\int_{\R}Z_1^T D\p_x(w_x\Omega N)\,dx,\label{eq4.15}\\
J_2&=\int_{\R}Z_3^T D\p_x^2 N\,dx+\int_{\R}w_x Z_2^T D'\p_x^2 N\,dx-\int_{\R}Z_2^T D\p_x^2(w_x\Omega N)\,dx.\label{eq4.16}
\end{align}
By Taylor expansion around $Z=0$,
\begin{align}\label{eq4.17}
N(w,Z)=\frac{1}{2}D_Z^2 N(w,0)[Z,Z]+\widetilde{N}(w,Z),\qquad |\widetilde{N}(w,Z)|\le C|Z|^3,
\end{align}
and the derivatives of $N$ and $w_x\Omega N$ satisfy
\begin{align}
|N|+|\Omega N|&\le C|Z|^2,\qquad |\p_x N|\le C\bigl(|Z||Z_1|+w_x|Z|^2\bigr),\label{eq4.18}\\
|\p_x^2 N|&\le C\bigl(|Z||Z_2|+|Z_1|^2+w_x|Z||Z_1|+(|w_{xx}|+w_x^2)|Z|^2\bigr),\nonumber\\
|\p_x(w_x\Omega N)|&\le C\bigl(w_x|Z||Z_1|+(|w_{xx}|+w_x^2)|Z|^2\bigr),\label{eq4.19}\\
|\p_x^2(w_x\Omega N)|&\le C\bigl[w_x(|Z||Z_2|+|Z_1|^2)+(|w_{xx}|+w_x^2)|Z||Z_1|+(|w_{xxx}|+w_x|w_{xx}|+w_x^3)|Z|^2\bigr].\nonumber
\end{align}
Substituting \eqref{eq4.17} into the first integral of $J_0$ in \eqref{eq4.14} produces the quadratic part $T_3(t)$ defined in \eqref{eq4.3}--\eqref{eq4.4}, while the cubic remainder and the two $w_x$-terms in $J_0$ satisfy
\[
\Bigl|\int_{\R}Z_1^T D\widetilde{N}\,dx\Bigr|+\Bigl|\int_{\R}w_x Z_0^T(D'-D\Omega)N\,dx\Bigr|\le C\varepsilon_b^2\mathcal{D}_2(t)+C\varepsilon_b\mathcal{G}_2(t).
\]
Similarly, \eqref{eq4.18}--\eqref{eq4.19} and $\|Z\|_{L^\infty}+\|Z_1\|_{L^\infty}\le C\varepsilon_b$ yield $|J_1|+|J_2|\le C\varepsilon_b(\mathcal{D}_2(t)+\mathcal{G}_2(t))$. Choosing $\eta>0$ small and then $\delta+\varepsilon_b$ small in \eqref{eq4.11}--\eqref{eq4.16} gives \eqref{eq4.2}, and integrating over $[0,\tau]$ with \eqref{eq2.5} yields \eqref{eq4.5}.
\end{proof}

\section{Quartic wave-interaction functional and \texorpdfstring{$H^2$}{H2} closure}\label{sec5}

\begin{Lem}\label{lem5.1}
Assume the bootstrap bound \eqref{eq4.1}. For $0\le\tau\le T$, define
\[
X_\tau:=\int_0^\tau\bigl(\mathcal{D}_2(t)+\mathcal{G}_2(t)\bigr)\,dt,\qquad Y_\tau:=\sum_{i<j}\int_0^\tau\int_{\R}z_i(t,x)^2 z_j(t,x)^2\,dx\,dt.
\]
If $\delta+\varepsilon_b$ is sufficiently small, then
\begin{align}\label{eq5.1}
Y_\tau\le C\mathcal{E}_2(0)^2+C\varepsilon_b^2 X_\tau+C\varepsilon_b^2\delta^3,
\end{align}
and for every $\eta>0$, the cubic flux integral in \eqref{eq4.3} satisfies
\begin{align}\label{eq5.2}
\Bigl|\int_0^\tau T_3(t)\,dt\Bigr|\le(\eta+C\varepsilon_b)X_\tau+C_\eta\mathcal{E}_2(0)^2+C_\eta\varepsilon_b^2\delta^3.
\end{align}
Consequently, for $\delta+\varepsilon_b$ sufficiently small, there exist constants $c,C_0>0$ independent of $T$ and $\tau$ such that
\begin{align}\label{eq5.3}
\sup_{0\le t\le\tau}\mathcal{E}_2(t)+cX_\tau\le C_0\bigl(\mathcal{E}_2(0)+\delta^3\bigr).
\end{align}
\end{Lem}

\begin{proof}
Write the $r$-th component of \eqref{eq4.7} as $z_{r,t}+\lambda_r(w)z_{r,x}-z_{r,xx}=S_r$, where
\begin{align}\label{eq5.4}
S=-w_x CZ+2w_x\Omega Z_x+w_x^2\Theta Z-L\p_x N_0+w_x^2 LU_c''(w).
\end{align}
For each pair $i<j$, define the one-sided quadratic tails and the directed interaction functional
\[
P_j^-(x):=\int_{-\infty}^x z_j(y)^2\,dy,\quad P_i^+(x):=\int_x^\infty z_i(y)^2\,dy,\quad \mathcal{I}_{ij}(t):=\iint_{x>y}z_i(x)^2 z_j(y)^2\,dx\,dy.
\]
Testing the $z_i$ equation against $2z_i P_j^-$, testing the $z_j$ equation against $2z_j P_i^+$, summing the two relations, and integrating by parts using $(P_j^-)_x=z_j^2$ and $(P_i^+)_x=-z_i^2$, we obtain the exact identity
\begin{align}\label{eq5.5}
&\mathcal{I}_{ij}'(t)+\int_{\R}(\lambda_j-\lambda_i)z_i^2 z_j^2\,dx+2\int_{\R}z_{i,x}^2 P_j^-\,dx+2\int_{\R}z_{j,x}^2 P_i^+\,dx\nonumber\\
&=\int_{\R}\bigl[z_i^2(z_j^2)_x-(z_i^2)_x z_j^2\bigr]\,dx+\iint_{x>y}\bigl(\lambda_{i,x}(x)+\lambda_{j,y}(y)\bigr)z_i(x)^2 z_j(y)^2\,dx\,dy\nonumber\\
&\quad+2\int_{\R}z_i S_i P_j^-\,dx+2\int_{\R}z_j S_j P_i^+\,dx.
\end{align}
By \eqref{eq4.1} and Sobolev embedding, the tails and perturbations obey
\begin{align}\label{eq5.6}
0\le P_j^-, P_i^+\le C\varepsilon_b^2,\qquad \|Z\|_{L^\infty}+\|Z_x\|_{L^\infty}\le C\varepsilon_b.
\end{align}
By Young's inequality, the viscous boundary integral on the right-hand side of \eqref{eq5.5} satisfies
\begin{align}\label{eq5.7}
\Bigl|\int_{\R}\bigl[z_i^2(z_j^2)_x-(z_i^2)_x z_j^2\bigr]\,dx\Bigr|\le\frac{\theta_0}{8}\int_{\R}z_i^2 z_j^2\,dx+C\bigl(\|z_i\|_{L^\infty}^2\|z_{j,x}\|_{L^2}^2+\|z_j\|_{L^\infty}^2\|z_{i,x}\|_{L^2}^2\bigr),
\end{align}
whose spacetime integral is bounded by $\frac{\theta_0}{8}Y_\tau+C\varepsilon_b^2 X_\tau$. Since $\lambda_{r,x}=w_x\lambda_r'(w)$, the variable-speed double integral in \eqref{eq5.5} is bounded by $C\varepsilon_b^2\mathcal{G}_2(t)$, and hence by $C\varepsilon_b^2 X_\tau$ after time integration. Likewise, using \eqref{eq5.6}, Cauchy--Schwarz, \eqref{eq2.3}--\eqref{eq2.5}, and $\sup_{t\ge 0}\|w_x(t)\|_{L^\infty}^{1/2}\le C\delta$, the four linear and profile source terms in $S$ contribute at most $C\varepsilon_b^2 X_\tau+C\varepsilon_b^2\delta^3$ to the time integral of \eqref{eq5.5}.

To bound the nonlinear flux terms in \eqref{eq5.5}, note first that for $i\neq j$ and any component $b$, integration by parts and pairing $z_i z_j$ against $z_b^2$ give
\begin{align}\label{eq5.8}
\int_{\R}z_b^4\,dx\le C\varepsilon_b^3\|z_{b,x}\|_{L^2},\qquad \int_{\R}z_j^2|z_i|z_b^2\,dx\le C\varepsilon_b^2\|z_{b,x}\|_{L^2}\Bigl(\int_{\R}z_i^2 z_j^2\,dx\Bigr)^{1/2}.
\end{align}
\noindent Indeed,
\begin{equation}
    \int_{\mathbb{R}} z_j^2\left|z_i\right| z_b^2 d x \leq\left\|z_i z_j\right\|_{L^2}\left\|z_j z_b^2\right\|_{L^2} \leq\left\|z_i z_j\right\|_{L^2}\left\|z_j\right\|_{L^2}\left\|z_b\right\|_{L^{\infty}}^2 \leq C \varepsilon_b^2\left\|z_{b, x}\right\|_{L^2}\left\|z_i z_j\right\|_{L^2}.
\end{equation}
\noindent Thus, the time integral of $\int_\R z_j^2|z_i|z_b^2\,dx$ over $[0,\tau]$ is bounded by $C\varepsilon_b^2 X_\tau^{1/2}Y_\tau^{1/2}$. Next, for any smooth function $a(w)$ with bounded derivative on $[w_-,w_+]$, we claim that
\begin{align}\label{eq5.9}
\sum_{b,c=1}^n\int_0^\tau\Bigl(\Bigl|\int_{\R}P_j^- a(w)z_{i,x}z_b z_c\,dx\Bigr|+\Bigl|\int_{\R}P_i^+ a(w)z_{j,x}z_b z_c\,dx\Bigr|\Bigr)\,dt\le C\varepsilon_b^2 X_\tau^{1/2}Y_\tau^{1/2}+C\varepsilon_b^3 X_\tau.
\end{align}
Indeed, if $b\neq c$, Cauchy--Schwarz in spacetime bounds the integral directly by $C\varepsilon_b^2 X_\tau^{1/2}Y_\tau^{1/2}$. If $b=c=r\neq i$ in the lower-tail integral, integrating by parts once gives
\[
\int_{\R}P_j^- a z_{i,x}z_r^2\,dx=-2\int_{\R}P_j^- a z_i z_r z_{r,x}\,dx-\int_{\R}\bigl(z_j^2 a+P_j^- w_x a'\bigr)z_i z_r^2\,dx,
\]
where the first term has $r\neq i$, the $z_j^2 a z_i z_r^2$ term is controlled by \eqref{eq5.8}, and the $w_x a'$ term is bounded by $C\varepsilon_b^3\mathcal{G}_2(t)$. If $b=c=i$ in the lower-tail integral, integrating by parts gives
\[
\int_{\R}P_j^- a z_{i,x}z_i^2\,dx=-\frac{1}{3}\int_{\R}\bigl(z_j^2 a+P_j^- w_x a'\bigr)z_i^3\,dx,
\]
which is again controlled by \eqref{eq5.8} and $C\varepsilon_b^3\mathcal{G}_2(t)$. The upper-tail integral with $P_i^+$ is treated identically, proving \eqref{eq5.9}.

Now, using \eqref{eq4.8} and integrating by parts, we write the nonlinear flux pairings in \eqref{eq5.5} as
\begin{align}\label{eq5.10}
2\int_{\R}P_j^- z_i(-L\p_x N_0)_i\,dx&=2\int_{\R}z_j^2 z_i N_i\,dx+2\int_{\R}P_j^- z_{i,x}N_i\,dx-2\int_{\R}P_j^- z_i w_x(\Omega N)_i\,dx,\nonumber\\
2\int_{\R}P_i^+ z_j(-L\p_x N_0)_j\,dx&=-2\int_{\R}z_i^2 z_j N_j\,dx+2\int_{\R}P_i^+ z_{j,x}N_j\,dx-2\int_{\R}P_i^+ z_j w_x(\Omega N)_j\,dx.
\end{align}
Decomposing $N=\frac{1}{2}D_Z^2 N(w,0)[Z,Z]+\widetilde{N}$ as in \eqref{eq4.17} and applying \eqref{eq5.8}--\eqref{eq5.9} to the quadratic part and Gagliardo--Nirenberg to the cubic remainder $\widetilde{N}$, we see that the time integral of \eqref{eq5.10} is bounded by $C\varepsilon_b^2 X_\tau^{1/2}Y_\tau^{1/2}+C\varepsilon_b^3 X_\tau$. Summing \eqref{eq5.5} over $i<j$, integrating over $[0,\tau]$ with $\lambda_j-\lambda_i\ge\theta_0$ and $\sum_{i<j}\mathcal{I}_{ij}(0)\le C\mathcal{E}_2(0)^2$, and applying Young's inequality yields \eqref{eq5.1}.

Finally, to bound the cubic flux integral in \eqref{eq4.3}, consider a single term with $a(w)=T_{rs\ell}(w)$. If $r=s=\ell$, integration by parts gives
\[
\int_{\R}a(w)z_{r,x}z_r^2\,dx=-\frac{1}{3}\int_{\R}w_x a'(w)z_r^3\,dx,
\]
whose time integral is bounded by $C\varepsilon_b X_\tau$. If $s\neq\ell$, Cauchy--Schwarz gives
\begin{align}\label{eq5.11}
\int_0^\tau\int_{\R}|z_{r,x}z_s z_\ell|\,dx\,dt\le X_\tau^{1/2}Y_\tau^{1/2}.
\end{align}
If $s=\ell\neq r$, integrating by parts once gives
\[
\int_{\R}a(w)z_{r,x}z_s^2\,dx=-2\int_{\R}a(w)z_r z_s z_{s,x}\,dx-\int_{\R}w_x a'(w)z_r z_s^2\,dx,
\]
where the first term has $r\neq s$ and is bounded by \eqref{eq5.11}, and the second term is bounded by $C\varepsilon_b\mathcal{G}_2(t)$. Summing over $(r,s,\ell)$ gives
\begin{align}\label{eq5.12}
\Bigl|\int_0^\tau T_3(t)\,dt\Bigr|\le C\varepsilon_b X_\tau+CX_\tau^{1/2}Y_\tau^{1/2}.
\end{align}
Substituting \eqref{eq5.1} into \eqref{eq5.12} and applying Young's inequality proves \eqref{eq5.2}. Combining \eqref{eq4.5} and \eqref{eq5.2} for small $\eta$ and $\varepsilon_b$ then yields \eqref{eq5.3}. In particular, choosing $\delta_0$ and $\varepsilon_0$ small enough that
\begin{align}\label{eq5.13}
C_0\bigl(\mathcal{E}_2(0)+\delta^3\bigr)\le\frac{c_E}{4}\varepsilon_b^2,\qquad C_1\bigl(\varepsilon_0+\delta_0^{3/2}\bigr)\le\frac{\kappa}{4}
\end{align}
for the Sobolev embedding constant $C_1>0$ ensures that \eqref{eq5.3} closes the bootstrap bound \eqref{eq4.1} and keeps $U(t,x)$ inside $\mathcal{K}$.
\end{proof}

\section{Weighted \texorpdfstring{$L^p$}{Lp} energy estimates}\label{sec6}
In Sections~\ref{sec6}--\ref{sec9}, we fix an integer $p\ge 4$ and allow constants to depend on $p$. By the closed $H^2$ estimate of Section~\ref{sec5}, the global solution satisfies
\begin{align}\label{eq6.1}
\sup_{t\ge 0}\|Z(t)\|_{H^2}+\sup_{t\ge 0}\|Z(t)\|_{L^\infty}\le CA.
\end{align}
For $i\neq k$, we define the $p$-adapted characteristic weight $d_i(w)$ along $[w_-,w_+]$ by solving
\begin{align}\label{eq6.2}
\frac{d_i'(w)}{d_i(w)}=p\Omega_{ii}(w)+\frac{(p-1)\lambda_i'(w)-\Gamma_p}{\lambda_i(w)-w},\qquad d_i(w_-)=1,
\end{align}
and set $d_k(w)\equiv 1$. We define the weighted $L^p$ energy, diffusion, rarefaction damping, and cross-family interaction functionals by
\begin{align}\label{eq6.3}
\begin{aligned}
\mathcal{E}_p(t)&:=\frac{1}{p}\sum_{i=1}^n\int_{\R}d_i(w)|z_i|^p\,dx,&\qquad \mathcal{D}_p(t)&:=\sum_{i=1}^n\int_{\R}|z_i|^{p-2}z_{i,x}^2\,dx,\\
\mathcal{G}_p(t)&:=\sum_{i=1}^n\int_{\R}w_x|z_i|^p\,dx,&\qquad \mathcal{Y}_p(t)&:=\sum_{i<j}\int_{\R}\bigl(|z_i|^p z_j^2+z_i^2|z_j|^p\bigr)\,dx.
\end{aligned}
\end{align}
We also collect the $L^1(0,\infty)$ quantities from the $H^2$ closure into the integrable function
\begin{align}\label{eq6.4}
a_2(t):=\mathcal{D}_2(t)+\mathcal{G}_2(t)+\mathcal{Y}_2(t)+\int_{\R}w_x(t,x)^3\,dx,\qquad \mathcal{Y}_2(t):=\sum_{i<j}\int_{\R}z_i^2 z_j^2\,dx,\qquad \int_0^\infty a_2(t)\,dt\le CA^2.
\end{align}

\begin{Lem}\label{lem6.1}
For every fixed integer $p\ge 4$, one can choose $\Gamma_p>0$ large and then $\tilde{\delta}_p>0$ small such that for $0<\delta\le\tilde{\delta}_p$,
\begin{align}\label{eq6.5}
\frac{1}{2}\le d_i(w)\le 2,\qquad |d_i'(w)|+|d_i''(w)|\le C_p,\qquad w\in[w_-,w_+],\quad 1\le i\le n.
\end{align}
Moreover, if $A+\delta$ is sufficiently small, then
\begin{align}\label{eq6.6}
\mathcal{E}_p'(t)+c_p\bigl(\mathcal{D}_p(t)+\mathcal{G}_p(t)\bigr)\le C_p\mathcal{Y}_p(t)+C_p\int_{\R}w_x(t,x)^{p+1}\,dx.
\end{align}
\end{Lem}

\begin{proof}
Set $\pi_i:=pd_i C_{ii}+wd_i'-(d_i\lambda_i)'$. Since $C_{ii}=\Omega_{ii}(\lambda_i-w)+\lambda_i'$, \eqref{eq6.2} gives $\pi_k=p-1$ and $\pi_i=\Gamma_p d_i$ for $i\neq k$. By Young's inequality, for $\Gamma_p>0$ chosen sufficiently large and $\tilde{\delta}_p>0$ small enough that \eqref{eq6.5} holds, we obtain the pointwise coercivity bound
\begin{align}\label{eq6.7}
\frac{1}{p}\sum_{i=1}^n\pi_i|z_i|^p-\sum_{i\neq j}d_i|C_{ij}||z_i|^{p-1}|z_j|\ge c_p\sum_{i=1}^n|z_i|^p.
\end{align}

Rewriting the $i$-th component of \eqref{eq4.7} as $z_{i,t}+\lambda_i z_{i,x}-z_{i,xx}+w_x C_{ii}z_i=\widetilde{S}_i$, where
\begin{align}\label{eq6.8}
\widetilde{S}_i:=-w_x\sum_{\ell\neq i}C_{i\ell}z_\ell+2w_x(\Omega Z_x)_i+w_x^2(\Theta Z)_i-\p_x N_i-w_x(\Omega N)_i+w_x^2(LU_c'')_i,
\end{align}
multiplying by $d_i(w)|z_i|^{p-2}z_i$, summing over $i$, and integrating by parts over $\R$ yields
\begin{align}\label{eq6.9}
&\mathcal{E}_p'(t)+(p-1)\sum_{i=1}^n\int_{\R}d_i|z_i|^{p-2}z_{i,x}^2\,dx+\frac{1}{p}\sum_{i=1}^n\int_{\R}w_x\pi_i|z_i|^p\,dx\nonumber\\
&=\sum_{i=1}^n\int_{\R}d_i|z_i|^{p-2}z_i\widetilde{S}_i\,dx-2\sum_{i=1}^n\int_{\R}w_x d_i'|z_i|^{p-2}z_i z_{i,x}\,dx-\frac{1}{p}\sum_{i=1}^n\int_{\R}w_x^2 d_i''|z_i|^p\,dx.
\end{align}
Moving the off-diagonal drift terms to the left-hand side of \eqref{eq6.9} and applying \eqref{eq6.7} gives the positive damping $c_p\mathcal{G}_p(t)$. For the multiplier derivative terms on the right-hand side of \eqref{eq6.9}, Young's inequality gives, for any $\eta>0$,
\[
2\sum_{i=1}^n\Bigl|\int_{\R}w_x d_i'|z_i|^{p-2}z_i z_{i,x}\,dx\Bigr|+\frac{1}{p}\sum_{i=1}^n\Bigl|\int_{\R}w_x^2 d_i''|z_i|^p\,dx\Bigr|\le\eta\mathcal{D}_p(t)+C_p\eta^{-1}\|w_x\|_{L^\infty}\mathcal{G}_p(t).
\]
For the frame-derivative term $2w_x(\Omega Z_x)_i$, when $\ell\neq i$ we integrate by parts with $b_{i\ell}:=2d_i\Omega_{i\ell}$:
\begin{align}\label{eq6.10}
\int_{\R}w_x b_{i\ell}(w)|z_i|^{p-2}z_i z_{\ell,x}\,dx=-(p-1)\int_{\R}w_x b_{i\ell}|z_i|^{p-2}z_{i,x}z_\ell\,dx-\int_{\R}\bigl(w_{xx}b_{i\ell}+w_x^2 b_{i\ell}'\bigr)|z_i|^{p-2}z_i z_\ell\,dx.
\end{align}
Combining \eqref{eq6.10} with $|w_{xx}|\le C\delta w_x$ and Young's inequality, we obtain for the linear and profile source terms in \eqref{eq6.8} that
\begin{align*}
\sum_{i=1}^n\Bigl|\int_{\R}d_i|z_i|^{p-2}z_i\bigl[2w_x(\Omega Z_x)_i+w_x^2(\Theta Z)_i+w_x^2(LU_c'')_i\bigr]\,dx\Bigr|&\le\eta\mathcal{D}_p(t)+\bigl[\eta+C_p\eta^{-1}(\delta+\|w_x\|_{L^\infty})\bigr]\mathcal{G}_p(t)\\
&\quad+C_p\eta^{-(p-1)}\int_{\R}w_x^{p+1}\,dx.
\end{align*}

It remains to estimate the nonlinear flux terms $-\p_x N_i-w_x(\Omega N)_i$ in \eqref{eq6.8}. We split $N_i(w,Z)$ into its single-mode part along the $i$-th axis and its transverse mixed part:
\begin{align}\label{eq6.11}
f_i(w,s):=N_i(w,se_i),\qquad R_i(w,Z):=N_i(w,Z)-f_i(w,z_i).
\end{align}
Since $N(w,0)=D_Z N(w,0)=0$, the single-mode primitive
\[
V_{i,p}(w,s):=\int_0^s|\zeta|^{p-2}\zeta\p_\zeta f_i(w,\zeta)\,d\zeta
\]
satisfies
\begin{align}\label{eq6.12}
|V_{i,p}(w,s)|+|\p_w V_{i,p}(w,s)|\le C_p|s|^{p+1},
\end{align}
and the chain rule gives
\begin{align}\label{eq6.13}
|z_i|^{p-2}z_i\p_x\bigl(f_i(w,z_i)\bigr)=\p_x\bigl(V_{i,p}(w,z_i)\bigr)+w_x\bigl(|z_i|^{p-2}z_i\p_w f_i(w,z_i)-\p_w V_{i,p}(w,z_i)\bigr).
\end{align}
Integrating by parts in $x$ and using \eqref{eq6.1} and \eqref{eq6.12} therefore gives
\[
\Bigl|-\int_{\R}d_i|z_i|^{p-2}z_i\p_x\bigl(f_i(w,z_i)\bigr)\,dx\Bigr|\le C_p\int_{\R}w_x|z_i|^{p+1}\,dx\le C_p A\mathcal{G}_p(t).
\]
For the mixed part $R_i(w,Z)$, Taylor's theorem gives the factorization
\begin{align}\label{eq6.14}
R_i(w,Z)=\sum_{a\neq i}\sum_{b=1}^n z_a z_b M_{iab}(w,Z),\qquad |M_{iab}|+|\p_w M_{iab}|+|D_Z M_{iab}|\le C.
\end{align}
Integrating the outer $\p_x$ by parts gives
\begin{align}\label{eq6.15}
-\int_{\R}d_i|z_i|^{p-2}z_i\p_x R_i\,dx=\int_{\R}w_x d_i'|z_i|^{p-2}z_i R_i\,dx+(p-1)\int_{\R}d_i|z_i|^{p-2}z_{i,x}R_i\,dx.
\end{align}
Since $|R_i(w,Z)|\le CA\sum_{a\neq i}|z_a|$, the first integral on the right-hand side of \eqref{eq6.15} and the lower-order flux term satisfy
\[
\Bigl|\int_{\R}w_x d_i'|z_i|^{p-2}z_i R_i\,dx\Bigr|+\Bigl|\int_{\R}d_i w_x|z_i|^{p-2}z_i(\Omega N)_i\,dx\Bigr|\le C_p A\mathcal{G}_p(t).
\]
For the second integral in \eqref{eq6.15}, Cauchy--Schwarz, Young's inequality, and the elementary bound for $p\ge 4$,
\begin{align}\label{eq6.16}
|z_i|^{p-2}|z_a|^4\le |z_i|^p|z_a|^2+|z_i|^2|z_a|^p\qquad (a\neq i),
\end{align}
yield
\[
\Bigl|\int_{\R}d_i|z_i|^{p-2}z_{i,x}R_i\,dx\Bigr|\le\eta\int_{\R}|z_i|^{p-2}z_{i,x}^2\,dx+C_p\eta^{-1}\sum_{a\neq i}\sum_{b=1}^n\int_{\R}|z_i|^{p-2}z_a^2 z_b^2\,dx\le\eta\mathcal{D}_p(t)+C_p\eta^{-1}\mathcal{Y}_p(t).
\]
Choosing $\eta>0$ small and then $A+\delta$ small in \eqref{eq6.9} proves \eqref{eq6.6}.
\end{proof}

\section{Mixed-power wave-interaction estimates}\label{sec7}
For the fixed integer $p\ge 4$, $i<j$, and $(r,s)\in\{(p,2),(2,p)\}$, we define the directed mixed-power interaction functionals and their one-sided tails by
\begin{align}\label{eq7.1}
\mathcal{I}_{ij}^{r,s}(t):=\iint_{x>y}|z_i(t,x)|^r|z_j(t,y)|^s\,dx\,dy,\quad P_{j,s}^-(x):=\int_{-\infty}^x|z_j(y)|^s\,dy,\quad P_{i,r}^+(x):=\int_x^\infty|z_i(y)|^r\,dy,
\end{align}
so that
\[
\mathcal{I}_{ij}^{r,s}(t)=\int_{\R}|z_i|^r P_{j,s}^-\,dx=\int_{\R}|z_j|^s P_{i,r}^+\,dx.
\]
Summing over both orientations $(p,2)$ and $(2,p)$, we set
\[
\mathcal{I}_p(t):=\sum_{i<j}\bigl(\mathcal{I}_{ij}^{p,2}(t)+\mathcal{I}_{ij}^{2,p}(t)\bigr),
\]
and denote the associated nonnegative cumulative diffusion integral by
\begin{align}\label{eq7.2}
\mathcal{B}_p(t):=\sum_{i<j}\int_{\R}\Bigl(p(p-1)P_{j,2}^-|z_i|^{p-2}z_{i,x}^2+2P_{i,p}^+ z_{j,x}^2+2P_{j,p}^- z_{i,x}^2+p(p-1)P_{i,2}^+|z_j|^{p-2}z_{j,x}^2\Bigr)\,dx.
\end{align}

\begin{Lem}\label{lem7.1}
Assume that the global $H^2$ bound of Lemma~\ref{lem5.1} holds. For every $\eta>0$, there exists $0<\tilde{a}_p\le\tilde{\delta}_p$ such that whenever $A+\delta\le\tilde{a}_p$,
\begin{align}\label{eq7.3}
\mathcal{I}_p'(t)+\frac{\theta_0}{2}\mathcal{Y}_p(t)+\mathcal{B}_p(t)\le\eta\bigl(\mathcal{D}_p(t)+\mathcal{G}_p(t)\bigr)+C_p\eta^{-1}a_2(t)\mathcal{E}_p(t)+C_p\eta^{-1}\int_{\R}w_x^{p+1}\,dx.
\end{align}
\end{Lem}

\begin{proof}
Differentiating \eqref{eq7.1} along $z_{\ell,t}+\lambda_\ell z_{\ell,x}-z_{\ell,xx}=S_\ell$ and integrating by parts over $x>y$ gives
\begin{align}\label{eq7.4}
&(\mathcal{I}_{ij}^{r,s})'+\int_{\R}(\lambda_j-\lambda_i)|z_i|^r|z_j|^s\,dx+r(r-1)\int_{\R}P_{j,s}^-|z_i|^{r-2}z_{i,x}^2\,dx+s(s-1)\int_{\R}P_{i,r}^+|z_j|^{s-2}z_{j,x}^2\,dx\nonumber\\
&=\int_{\R}\bigl[|z_i|^r(|z_j|^s)_x-(|z_i|^r)_x|z_j|^s\bigr]\,dx+\int_{\R}\lambda_{i,x}|z_i|^r P_{j,s}^-\,dx+\int_{\R}\lambda_{j,x}|z_j|^s P_{i,r}^+\,dx\nonumber\\
&\quad+r\int_{\R}P_{j,s}^-|z_i|^{r-2}z_i S_i\,dx+s\int_{\R}P_{i,r}^+|z_j|^{s-2}z_j S_j\,dx.
\end{align}
Integrating the first term on the right-hand side of \eqref{eq7.4} by parts and applying Young's inequality with $\|z_\ell\|_{L^\infty}\le CA$ gives
\[
2p\sum_{i<j}\Bigl(\Bigl|\int_{\R}|z_i|^{p-2}z_i z_{i,x}z_j^2\,dx\Bigr|+\Bigl|\int_{\R}z_i^2|z_j|^{p-2}z_j z_{j,x}\,dx\Bigr|\Bigr)\le\frac{\theta_0}{8}\mathcal{Y}_p(t)+C_p A^2\mathcal{D}_p(t).
\]
The cumulative tails satisfy
\begin{align}\label{eq7.5}
\|P_{\ell,2}^\pm\|_{L^\infty}\le\|z_\ell\|_{L^2}^2\le CA^2,\qquad \|P_{\ell,p}^\pm\|_{L^\infty}\le\|z_\ell\|_{L^p}^p\le C_p\mathcal{E}_p(t).
\end{align}
Hence, since $|\lambda_{\ell,x}|\le Cw_x$, the speed-derivative integrals in \eqref{eq7.4} are bounded by $C_p A^2\mathcal{G}_p(t)+C_p\mathcal{E}_p(t)\mathcal{G}_2(t)$.

We now substitute $S_\ell$ from \eqref{eq5.4} into the last two integrals of \eqref{eq7.4}. Using \eqref{eq7.5}, $|z_i|^{p-1}|z_\ell|\le C_p(|z_i|^p+|z_\ell|^p)$, and $|z_i z_\ell|\le\frac{1}{2}(z_i^2+z_\ell^2)$, we bound the linear drift and Burgers residual terms on the $p$-side ($r=p$ or $s=p$) and on the quadratic side ($r=2$ or $s=2$) by
\begin{align*}
\Bigl|\int_{\R}P_{j,2}^-|z_i|^{p-2}z_i\bigl[-w_x(CZ)_i+w_x^2(\Theta Z)_i+w_x^2(LU_c'')_i\bigr]\,dx\Bigr|&\le\bigl(\eta+C_p A^2\bigr)\mathcal{G}_p(t)+C_p\eta^{-1}\int_{\R}w_x^{p+1}\,dx,\\
\Bigl|\int_{\R}P_{j,p}^- z_i\bigl[-w_x(CZ)_i+w_x^2(\Theta Z)_i+w_x^2(LU_c'')_i\bigr]\,dx\Bigr|&\le C_p\mathcal{E}_p(t)\Bigl(\mathcal{G}_2(t)+\int_{\R}w_x^3\,dx\Bigr),
\end{align*}
with identical bounds for the upper-tail integrals involving $P_{i,2}^+$ and $P_{i,p}^+$.

For the frame-derivative term $2w_x(\Omega Z_x)_i$, Cauchy--Schwarz on the quadratic side gives
\[
\Bigl|\int_{\R}P_{j,p}^- w_x z_i(\Omega Z_x)_i\,dx\Bigr|\le C_p\mathcal{E}_p(t)\mathcal{G}_2(t).
\]
On the $p$-side, when $\ell=i$ we apply Young's inequality directly, and when $\ell\neq i$ we integrate by parts in $x$ using $(P_{j,2}^-)_x=z_j^2$ ($j\neq i$) together with \eqref{eq6.16}, obtaining
\begin{align*}
\Bigl|\int_{\R}P_{j,2}^- w_x|z_i|^{p-2}z_i(\Omega Z_x)_i\,dx\Bigr|&\le\eta\mathcal{D}_p(t)+C_p\eta^{-1}A^2\bigl(\delta+\|w_x\|_{L^\infty}\bigr)\mathcal{G}_p(t)+C_p\|w_x\|_{L^\infty}\mathcal{Y}_p(t).
\end{align*}

It remains to estimate the nonlinear flux $-L\p_x N_0=-\p_x N-w_x\Omega N$ in $S_\ell$. On the lower $p$-side, we use the splitting $N_i=f_i(w,z_i)+R_i(w,Z)$ from \eqref{eq6.11}. By \eqref{eq6.12}--\eqref{eq6.13} and $(P_{j,2}^-)_x=z_j^2$, integrating by parts gives
\begin{align*}
\Bigl|-\int_{\R}P_{j,2}^-|z_i|^{p-2}z_i\p_x\bigl(f_i(w,z_i)\bigr)\,dx\Bigr|&=\Bigl|\int_{\R}z_j^2 V_{i,p}\,dx-\int_{\R}P_{j,2}^- w_x\bigl(|z_i|^{p-2}z_i\p_w f_i-\p_w V_{i,p}\bigr)\,dx\Bigr|\\
&\le C_p A\mathcal{Y}_p(t)+C_p A^3\mathcal{G}_p(t).
\end{align*}
For the mixed part $R_i$, integrating the outer $\p_x$ by parts and applying \eqref{eq6.14}, \eqref{eq6.16}, and Young's inequality yields
\begin{align*}
\Bigl|-\int_{\R}P_{j,2}^-|z_i|^{p-2}z_i\p_x R_i\,dx\Bigr|&=\Bigl|\int_{\R}z_j^2|z_i|^{p-2}z_i R_i\,dx+(p-1)\int_{\R}P_{j,2}^-|z_i|^{p-2}z_{i,x}R_i\,dx\Bigr|\\
&\le\eta\mathcal{D}_p(t)+C_p\eta^{-1}A\mathcal{Y}_p(t),
\end{align*}
while the $-w_x(\Omega N)_i$ part is bounded by $C_p A^3\mathcal{G}_p(t)$. The upper $p$-side with $P_{i,2}^+$ is estimated identically.

Finally, consider the lower quadratic side paired with $P_{j,p}^- z_i$, where $\|P_{j,p}^-\|_{L^\infty}\le C_p\mathcal{E}_p(t)$ and $j\neq i$. For the single-mode part $f_i$, using the quadratic primitive $V_{i,2}(w,s):=\int_0^s\zeta\p_\zeta f_i(w,\zeta)\,d\zeta$ and integrating by parts gives
\[
\Bigl|-\int_{\R}P_{j,p}^- z_i\p_x\bigl(f_i(w,z_i)\bigr)\,dx\Bigr|\le C_p\int_{\R}|z_j|^p|z_i|^3\,dx+C_p\mathcal{E}_p(t)\int_{\R}w_x|z_i|^3\,dx\le C_p A\mathcal{Y}_p(t)+C_p A\mathcal{E}_p(t)\mathcal{G}_2(t).
\]
For the mixed part $R_i=\sum_{a\neq i}\sum_{b=1}^n z_a z_b M_{iab}(w,Z)$, expanding $z_i\p_x(z_a z_b M_{iab})$ gives
\begin{align}\label{eq7.6}
z_i\p_x(z_a z_b M_{iab})=z_i z_{a,x}z_b M_{iab}+z_i z_a z_{b,x}M_{iab}+w_x z_i z_a z_b\p_w M_{iab}+z_i z_a z_b\bigl(D_Z M_{iab}\cdot Z_x\bigr).
\end{align}
Since $a\neq i$, the last three terms in \eqref{eq7.6} all contain the distinct product $z_i z_a$ or $w_x$:
\begin{align}\label{eq7.7}
&\Bigl|\int_{\R}P_{j,p}^- z_i z_a z_{b,x}M_{iab}\,dx\Bigr|+\Bigl|\int_{\R}P_{j,p}^- w_x z_i z_a z_b\p_w M_{iab}\,dx\Bigr|+\Bigl|\int_{\R}P_{j,p}^- z_i z_a z_b(D_Z M_{iab}\cdot Z_x)\,dx\Bigr|\nonumber\\
&\le C_p\mathcal{E}_p(t)\bigl(\mathcal{D}_2(t)+\mathcal{Y}_2(t)+A\mathcal{G}_2(t)\bigr).
\end{align}
In the first term of \eqref{eq7.6}, if $b\neq i$, it likewise contains the distinct product $z_i z_b$ and is bounded by the right-hand side of \eqref{eq7.7}. If $b=i$, we integrate by parts in $x$ before taking the supremum of $P_{j,p}^-$:
\begin{align*}
\int_{\R}P_{j,p}^- M_{iab}z_{a,x}z_i^2\,dx&=-2\int_{\R}P_{j,p}^- M_{iab}z_a z_i z_{i,x}\,dx-\int_{\R}|z_j|^p M_{iab}z_a z_i^2\,dx\\
&\quad-\int_{\R}P_{j,p}^-\bigl(w_x\p_w M_{iab}+D_Z M_{iab}\cdot Z_x\bigr)z_a z_i^2\,dx,
\end{align*}
where the first and third integrals on the right satisfy \eqref{eq7.7} and the second integral is bounded by $C_p A\mathcal{Y}_p(t)$. Combining these bounds with \eqref{eq7.4}, summing over $i<j$, and taking $A+\delta$ sufficiently small yields \eqref{eq7.3}.
\end{proof}

\section{Closure of the \texorpdfstring{$L^p$}{Lp} Lyapunov functional}\label{sec8}
By Fubini's theorem and \eqref{eq6.5}, the interaction functional $\mathcal{I}_p(t)$ satisfies
\begin{align}\label{eq8.1}
0\le\mathcal{I}_p(t)\le C_p\|Z(t)\|_{L^2}^2\mathcal{E}_p(t)\le C_p A^2\mathcal{E}_p(t).
\end{align}

\begin{Lem}\label{lem8.1}
For each fixed integer $p\ge 4$, there exist constants $K_p>0$ and $a_p>0$ such that for $A+\delta\le a_p$, the Lyapunov functional
\[
\mathcal{H}_p(t):=\mathcal{E}_p(t)+K_p\mathcal{I}_p(t)
\]
satisfies
\begin{align}
\mathcal{E}_p(t)&\le\mathcal{H}_p(t)\le 2\mathcal{E}_p(t),\label{eq8.2}\\
\mathcal{H}_p'(t)+c_p\bigl(\mathcal{D}_p(t)+\mathcal{G}_p(t)+\mathcal{Y}_p(t)+\mathcal{B}_p(t)\bigr)&\le C_p a_2(t)\mathcal{H}_p(t)+C_p\int_{\R}w_x^{p+1}\,dx,\label{eq8.3}
\end{align}
and consequently
\begin{align}\label{eq8.4}
\sup_{t\ge 0}\mathcal{H}_p(t)+c_p\int_0^\infty\bigl(\mathcal{D}_p(t)+\mathcal{G}_p(t)+\mathcal{Y}_p(t)+\mathcal{B}_p(t)\bigr)\,dt\le C_p A^p.
\end{align}
\end{Lem}

\begin{proof}
Adding \eqref{eq6.6} and $K_p$ times \eqref{eq7.3}, choosing $K_p>0$ large enough that $\frac{1}{2}K_p\theta_0$ absorbs $C_p\mathcal{Y}_p(t)$ on the right-hand side of \eqref{eq6.6}, and then choosing $\eta>0$ in Lemma~\ref{lem7.1} and $a_p>0$ small enough that $K_p\eta\le\frac{1}{2}c_p$ and $C_p K_p A^2\le 1$ in \eqref{eq8.1}, we obtain \eqref{eq8.2} and \eqref{eq8.3}. Moreover, by \eqref{eq2.6},
\[
\mathcal{H}_p(0)\le C_p\|Z(0)\|_{L^\infty}^{p-2}\|Z(0)\|_{L^2}^2\le C_p A^p,\qquad \int_0^\infty\int_{\R}w_x^{p+1}\,dx\,dt\le C_p\delta^{2p-1}\le C_p A^p.
\]
Applying Gronwall's inequality to \eqref{eq8.3} using \eqref{eq6.4} therefore yields \eqref{eq8.4} and \eqref{eq2.10}.
\end{proof}

\section{Nonlinear Nash inequality and \texorpdfstring{$L^p$}{Lp} decay rates}\label{sec9}

\begin{Lem}\label{lem9.1}
(1)For any fixed integer $p\ge 4$ and $u\in H^1(\R)\cap L^p(\R)$,
\begin{align}\label{eq9.1}
\Bigl(\int_{\R}|u|^p\,dx\Bigr)^{\frac{p+2}{p-2}}\le C_p\Bigl(\int_{\R}|u|^2\,dx\Bigr)^{\frac{2p}{p-2}}\int_{\R}|u|^{p-2}u_x^2\,dx.
\end{align}
(2)Under the assumption $\|Z(t)\|_{L^2}\le CA$, it holds that
\begin{align}\label{eq9.2}
\mathcal{D}_p(t)\ge c_p A^{-\frac{4p}{p-2}}\mathcal{H}_p(t)^{\frac{p+2}{p-2}},
\end{align}
and the solution satisfies
\[
\mathcal{H}_p(t)\le C_p A^p(1+t)^{-(p-2)/4},\qquad \|W(t)\|_{L^p}\le C_p A(1+t)^{-\frac{1}{2}(\frac{1}{2}-\frac{1}{p})}\qquad \forall\,t\ge 0.
\]
\end{Lem}

\begin{proof}
Applying the one-dimensional Gagliardo--Nirenberg inequality \cite{Nash1958} to $|u|^{p/2}$ gives \eqref{eq9.1}. Summing \eqref{eq9.1} over $u=z_i$ ($1\le i\le n$) and using $\|Z(t)\|_{L^2}\le CA$ together with \eqref{eq8.2} yields \eqref{eq9.2}. By \eqref{eq2.7}, $A\ge\delta^{3/2}$, $0<\delta\le 1$, and $\frac{p+2}{4}\le p$ for $p\ge 2$, the profile source in \eqref{eq8.3} satisfies
\begin{align}\label{eq9.3}
\begin{aligned}
\int_{\mathbb{R}} w_x^{p+1} d x \leq C_p \delta^{2 p+1}\left(1+\delta^2 t\right)^{-p} & \leq C_p \delta^{2 p+1}\left(1+\delta^2 t\right)^{-\frac{p+2}{4}} \leq C_p \delta^{2 p+1}\left(\delta^2(1+t)\right)^{-\frac{p+2}{4}} \\
& =C_p \delta^{\frac{3 p}{2}}(1+t)^{-\frac{p+2}{4}} \leq C_p A^p(1+t)^{-\frac{p+2}{4}}.
\end{aligned}
\end{align}
Substituting \eqref{eq9.2} and \eqref{eq9.3} into \eqref{eq8.3} gives
\[
\mathcal{H}_p'(t)+c_p A^{-\frac{4p}{p-2}}\mathcal{H}_p(t)^{\frac{p+2}{p-2}}\le C_p a_2(t)\mathcal{H}_p(t)+C_p A^p(1+t)^{-\frac{p+2}{4}}.
\]
Multiplying by the integrating factor $\exp(-C_p\int_0^t a_2(s)\,ds)\simeq 1$ (by \eqref{eq6.4}) and comparing with the barrier $Y(t):=M_p A^p(1+t)^{-(p-2)/4}$ for $M_p>0$ sufficiently large yields $\mathcal{H}_p(t)\le C_p A^p(1+t)^{-(p-2)/4}$. Taking the $p$-th root and using $\|W(t)\|_{L^p}^p\le C_p\mathcal{H}_p(t)$ gives
\begin{align}\label{eq9.4}
\|W(t)\|_{L^p}\le C_p A(1+t)^{-\frac{p-2}{4p}}=C_p A(1+t)^{-\frac{1}{2}(\frac{1}{2}-\frac{1}{p})}.
\end{align}
Finally, for $p=3$, H\"older interpolation $\|W(t)\|_{L^3}\le\|W(t)\|_{L^2}^{1/3}\|W(t)\|_{L^4}^{2/3}\le CA(1+t)^{-1/12}$ completes the proof of \eqref{eq2.11}.
\end{proof}

\section{Global continuation and large-time asymptotics}\label{sec10}

\begin{Prop}\label{prop10.1}
Assume $W_0=U_0-\bar{U}(0,\cdot)\in H^2(\R)$ and $\overline{U_0(\R)}\Subset\mathcal{O}$. Then \eqref{eq1.2} admits a unique maximal strong solution on $[0,T_{\max})$ satisfying, for every $T<T_{\max}$,
\begin{align}\label{eq10.1}
W\in C([0,T];H^2(\R))\cap L^2(0,T;H^3(\R)),\qquad W_t\in L^2(0,T;H^1(\R)).
\end{align}
If
\begin{align}\label{eq10.2}
\sup_{0\le t<T_{\max}}\|W(t)\|_{H^2}<\infty\qquad\text{and}\qquad \overline{\bigcup_{0\le t<T_{\max}}U(t,\R)}\Subset\mathcal{O},
\end{align}
then $T_{\max}=\infty$. Under the smallness assumptions of Theorem~\ref{thm2.1}, $T_{\max}=\infty$ and
\[
\lim_{t\to\infty}\|W(t)\|_{L^\infty}=0,\qquad \lim_{t\to\infty}\sup_{x\in\R}|U(t,x)-U^R(x/t)|=0.
\]
\end{Prop}

\begin{proof}
By a standard smooth cutoff of $F$ outside a compact neighborhood of $\overline{U_0(\R)}\cup K_r\Subset\mathcal{O}$ and fixed-point contraction for the semilinear parabolic system \eqref{eq4.6} in $H^2(\R)$, there exists a unique maximal strong solution satisfying \eqref{eq10.1} and the continuation criterion \eqref{eq10.2}. Defining the bootstrap exit time
\begin{align}\label{eq10.3}
T_*:=\sup\Bigl\{T<T_{\max}:\sup_{0\le t\le T}\|Z(t)\|_{H^2}\le\varepsilon_b,\ U([0,T]\times\R)\subset\mathcal{K}\Bigr\},
\end{align}
we see from \eqref{eq5.3} and \eqref{eq5.13} that on $[0,T_*]$ one has $\sup_{0\le t\le T_*}\|Z(t)\|_{H^2}\le\frac{1}{2}\varepsilon_b$ and $\operatorname{dist}(U(t,x),K_r)\le\|W(t)\|_{L^\infty}\le\frac{1}{2}\kappa$. By continuity and \eqref{eq10.2}, $T_*=T_{\max}=\infty$, which establishes global existence and \eqref{eq2.9}.

To verify \eqref{eq2.12}, set $G(t):=\|W_x(t)\|_{L^2}^2$. By \eqref{eq2.9}, $G\in L^1(0,\infty)$. Differentiating $G(t)$ and using \eqref{eq4.6} gives
\begin{align*}
G'(t)=-2\int_{\R}W_{xx}\cdot W_t\,dx&=-2\|W_{xx}\|_{L^2}^2+2\int_{\R}W_{xx}\cdot\bigl(A_F(\bar{U})W\bigr)_x\,dx\\
&\quad+2\int_{\R}W_{xx}\cdot\p_x N_0\,dx-2\int_{\R}W_{xx}\cdot U_c''(w)w_x^2\,dx,
\end{align*}
and hence, by Cauchy--Schwarz and $\|W\|_{H^2}\le C\varepsilon_b$,
\begin{align}\label{eq10.4}
|G'(t)|\le C\Bigl(\|W_x\|_{H^1}^2+\int_{\R}w_x|W|^2\,dx+\|w_x\|_{L^\infty}^2\|W\|_{L^2}^2+\int_{\R}w_x^3\,dx\Bigr).
\end{align}
By \eqref{eq2.5}, \eqref{eq2.7}, and \eqref{eq2.9}, the right-hand side of \eqref{eq10.4} belongs to $L^1(0,\infty)$, so $G\in W^{1,1}(0,\infty)$ and $\lim_{t\to\infty}\|W_x(t)\|_{L^2}=0$. By the Gagliardo--Nirenberg inequality,
\begin{align}\label{eq10.5}
\|W(t)\|_{L^\infty}^2\le 2\|W(t)\|_{L^2}\|W_x(t)\|_{L^2}\longrightarrow 0\qquad\text{as }t\to\infty.
\end{align}
Since $\sup_{x\in\R}|\bar{U}(t,x)-U^R(x/t)|\to 0$ as $t\to\infty$ by the classical Burgers comparison estimate \cite{IlinOleinik1960,Xin1988,HattoriNishihara1991}, combining this with \eqref{eq10.5} yields \eqref{eq2.12} and completes the proof of Theorem~\ref{thm2.1}.
\end{proof}

\bigskip
\noindent\textbf{Competing interests.} The authors declare that they have no competing interests.

\medskip
\noindent\textbf{Data availability.} Data sharing is not applicable to this article as no datasets were generated or analyzed during the current study.

\end{document}